\documentclass[11pt,a4paper]{article}
\usepackage{pb-diagram}
\makeatletter
\renewcommand\thesection{\@Roman\c@section}
\renewcommand\thesubsection{\thesection.\@arabic\c@subsection}
\makeatother

\newcommand{\sect}[1]{\setcounter{equation}{0}\section{#1}}

\newtheorem{theorem}{Theorem}[section]
\newtheorem{lemma}[theorem]{Lemma}
\newtheorem{corollary}[theorem]{Corollary}
\newtheorem{proposition}[theorem]{Proposition}
\newtheorem{definition}{Definition}[section]

\begin{document}
\title{\Large\bf Quantum Doubles From A Class Of Noncocommutative Weak Hopf
Algebras}
\author{ {\large Fang Li} \thanks{fangli@zju.edu.cn}\\
Department of Mathematics, Zhejiang
University\\Hangzhou, Zhejiang 310028, China\\
\\
{\large Yao-Zhong Zhang} \thanks{yzz@maths.uq.edu.au}\\
Department of Mathematics, University of Queensland\\
Brisbane, Qld 4072, Australia}
\maketitle

\begin{abstract}
The concept of biperfect (noncocommutative) weak Hopf algebras is
introduced and their properties are discussed.  A new type of
quasi-bicrossed products are constructed by means of weak Hopf
skew-pairs of the weak Hopf algebras which are generalizations of
the Hopf pairs introduced by Takeuchi. As a special case, the
quantum double of a finite dimensional biperfect
(noncocommutative) weak Hopf algebra is built.
Examples of quantum doubles from a Clifford monoid as
well as a noncommutative and noncocommutative weak Hopf algebra
are given, generalizing quantum doubles from a group and a
noncommutative and noncocommutative Hopf algebra, respectively.
Moreover, some characterisations of quantum doubles of finite
dimensional biperfect weak Hopf algebras are obtained.
\end{abstract}

\sect{Introduction}

In a recent work \cite{L1}, quantum doubles of finite dimensional Hopf
algebras and finite groups are generalized by one of the authors
to those of certain finite dimensional weak Hopf algebras and finite
Clifford monoids so as to obtain singular solutions of the quantum
Yang-Baxter equation. However, the procedure in \cite{L1}
is not suitable for noncocommutative weak Hopf algebras. So, it is interesting
to construct quantum doubles of noncocommutative weak Hopf algebras.
The aim of this paper is to give one class of noncocommutative weak Hopf
algebras from which their quantum doubles can be obtained.

As is known \cite{K}, bicrossed product is a fundamental tool to
construct the quantum double of a Hopf algebra.
Quasi-bicrossed product plays a similar role in \cite{L1} for the construction
of quantum doubles of certain weak Hopf algebras (in particular, finite Clifford
monoids). In \cite{L4}, the concept of weak Hopf pairs was introduced as a
generalization of the Hopf pairs of Takeuchi \cite{T}.
Using the weak Hopf skew-pairs, one type of quasi-bicrossed
products, which lie between general quasi-bicrossed products
and quantum quasi-doubles, were constructed when one of the two weak Hopf
algebras in the product is cocommutative. In section \ref{sec2} of this paper,
we will generalize the results and construct quasi-bicrossed products of
the weak Hopf skew-pairs corresponding to the case where both weak
Hopf algebras in the product are non-cocommutative (see Theorem
\ref{th2.6} below).

A bialgebra $H$ over a field $k$ is called a {\em weak Hopf algebra} \cite{L1}
if there exists $T\in Hom_{k}(H,H)$ such that $id\ast T\ast id=id$ and $T\ast
id\ast T=T$ where $\ast$ is the covolution product in $Hom_{k}(H,H)$; $T$
is called a {\em weak antipode} of $H$. Weak Hopf algebras lie between left
(resp. right) Hopf algebras and bialgebras. So far, two types of
such weak Hopf algebras have been found, which are the monoid algebra $kS$ of
 a regular monoid $S$ \cite{L1} and the almost quantum algebra $wsl_{q}(2)$
\cite{L3} (see also \cite{Isaac} for weak Hopf algebras corresponding to
$U_q[sl_n]$).

An application of weak Hopf algebras was found in the construction of
non-invertible solutions of the (quantum) Yang-Baxter equation
in \cite{L1,L3}. It was found
that for a finite dimensional cocommutative perfect weak Hopf algebra $H$
with an invertible weak antipode $T$, the quasi-bicrossed product $H^{op\ast}
\infty H$ (which is called the {\em quantum double} of $H$, denoted by $D(H)$),
is a quasi-braided almost bialgebra equipped with the quasi-R-matrix
$R=\sum^{n}_{i=1}(1\infty e_{i})\otimes(e^{i}\infty 1)\in D(H)\otimes D(H)$
where $\{e_{i}\}^{n}_{i=1}$ is a basis of $H$ as a vector space and
$\{e^{i}\}^{n}_{i=1}$ is its dual basis in $H^{op\ast}$. Then, $R$ is a
solution of the quantum Yang-Baxter equation. In \cite{L2}, it was shown
that this solution $R$ is von Neumann regular but not invertible in general.
An example of this solution was constructed
from the cocommutative perfect weak Hopf
algebra $H=kS$ for any finite Clifford monoid $S$.

Although the quantum double of a finite Clifford monoid is indeed a
generalization of the quantum double of a finite group,
the quantum doubles in \cite{L1} can not usually be regarded
as  generalizations of quantum doubles of Hopf algebras due to the
cocommutativity of the weak Hopf algebras considered in \cite{L1}.
The goal of this paper is to overcome this  so as to construct quantum doubles
of noncocommutative weak Hopf algebras.
 We will give one class of noncocommutative weak Hopf algebras from which their
quantum doubles can be obtained. Firstly, we introduce the concept of
biperfect weak Hopf algebras and discuss their properties.
Then we construct a new type of quasi-bicrossed products by means of
the weak Hopf skew-pairs of the weak Hopf algebras which are generalizations
of  the Hopf pairs introduced by Takeuchi \cite{L4}. As a special case,
the quantum double of a finite dimensional
biperfect (noncocommutative) weak Hopf algebra is built.
Examples of quantum doubles from a Clifford monoid and a noncommutative
and noncocommutative weak Hopf algebra
are given as generalizations of those from a group and a noncommutative and
noncocommutative Hopf algebra, respectively.
Moreover, we discuss some characterisations of quantum doubles of finite
dimensional biperfect weak Hopf algebras.

\sect{Preliminaries}
Throughout the paper, $k$ stands for a field. Some notations and definitions
unexplained here can be found in \cite{P}, \cite{K}, \cite{L1} and \cite{S}.
The word ``quantum quasi-double" of a weak Hopf algebra in \cite{L1}
will always be replaced with ``quantum double".

We recall \cite{L1} that a linear space $H$ is a {\em $k$-almost
bialgebra} if $(H,\mu,\eta)$ is a $k$-algebra and $(H,\Delta,\varepsilon)$
is a $k$-coalgebra with $\Delta (xy)=\Delta (x)\Delta (y)$ for $x$,
$y\in H$. If $K$ is a subalgebra and also a sub-coalgebra of $H$, then $K$
itself is an almost bialgebra, called as an {\em almost sub-bialgebra}
of $H$.

Combining formally with the definition of the weak Hopf algebras, we say in
\cite{L2} that an almost bialgebra $H$ is an {\em almost weak Hopf
algebra} if there exists $T\in Hom_{k}(H,H)$ such that $id\ast T\ast id=id$
and $T\ast id\ast T=T$, where $T$ is called an {\em almost weak
antipode} of $H$.

Let $H$ be an almost bialgebra. If there exists an $R\in H\otimes H$
such that for all $x\in H$, $\Delta^{op}(x)R=R\Delta(x)$, then $R$ is called
a  {\em universal quasi-R-matrix}; if simultaneously,
$(\Delta \otimes id)(R)=R_{13}R_{23}$ and $(id\otimes\Delta)(R)=R_{13}R_{12}$
are satisfied, then we call $H$ a {\em quasi-braided almost bialgebra}
with a {\em quasi-R-matrix} $R$ (see \cite{L1}). Moreover, if $H$ is a
bialgebra and $R$ is invertible, then $H$ is called a {\em braided bialgebra}.

Let $H$ be a bialgebra and $C$ a coalgebra. If $C$ is a left $H$-module
and $\Delta (hc)=\Delta (h)\Delta (c)$ for
every $h\in H$ and $c\in C$, then we call the coalgebra $C$ a {\em left
quasi-module-coalgebra} over $H$. Moreover,
if $\varepsilon(hc)=\varepsilon(h)\varepsilon(c)$, then $C$ is
called a {\em left module-coalgebra} over $H$. Right
quasi-module-coalgebra and right module-coalgebra can be defined similarly.

A pair $(X,A)$ of bialgebras over a field $k$ is called {\em
quasi-matched} (resp. {\em matched}) if there exist linear maps
$\alpha: A\otimes
X\longrightarrow X$ and $\beta: A\otimes X\longrightarrow A$ which turn $X$
into a left $A$-quasi-module-coalgebra (resp. a left
$A$-module-coalgebra) and $A$ into a right $X$-quasi-module-coalgebra (resp. a
right $X$-module-coalgebra), such that if one sets $\alpha (a\otimes
x)=a\triangleright x$, $\beta (a\otimes x)=a\triangleleft x$ then the following
conditions are satisfied:
\begin{equation}
a\triangleright (xy)=\sum_{(a)(x)}(a'\triangleright
x')((a''\triangleleft x'')\triangleright y),\label{e1}
\end{equation}
\begin{equation}
a\triangleright 1=\varepsilon (a)1,\label{e2}
\end{equation}
\begin{equation}
(ab)\triangleleft x=\sum_{(b)(x)}(a\triangleleft (b'\triangleright
x'))(b''\triangleleft x''),\label{e3}
\end{equation}
\begin{equation}
1\triangleleft x=\varepsilon (x)1,\label{e4}
\end{equation}
\begin{equation}
\sum_{(a)(x)}(a'\triangleleft x')\otimes (a''\triangleright
x'')=\sum_{(a)(x)}(a''\triangleleft x'')\otimes (a'\triangleright x'),\label{e5}
\end{equation}
for all $a$, $b\in A$ and $x$, $y\in X$, where $1$ is the identity of
$X$ and of $A$ respectively in (\ref{e2}) and (\ref{e4}).

For a quasi-matched (resp. matched) pair of bialgebras $(X,A)$, we know
from \cite{L1}  and \cite{K}  that there exists an almost bialgebra
(resp. a bialgebra) structure on the vector space $X\otimes A$ with
identity equal to $1\otimes 1$ such that its product is given by
\begin{equation}
(x\otimes a)(y\otimes b)=\sum_{(a)(y)}x(a'\triangleright y')\otimes
(a''\triangleleft y'')b,\label{e6}
\end{equation}
\\ its coproduct by
\begin{equation}
\Delta (x\otimes a)=\sum_{(a)(x)}(x'\otimes a')\otimes (x''\otimes a'')
\label{e7}
\end{equation}
\\ and its counit by
\begin{equation}
\varepsilon (x\otimes a)={\varepsilon}_{X}(x){\varepsilon}_{A}(a)\label{e8}
\end{equation}
for all $x$, $y\in X$, $a$, $b\in A$. Equipped with this almost
bialgebra (resp. bialgebra) structure, $X\otimes A$ is called the
{\em quasi-bicrossed product} (resp. {\em bicrossed product}) of $X$ and $A$,
and denoted as
$X\infty A$. Furthermore, the injective maps $i_{X}(x)=x\otimes 1$ and
$i_{A}(a)=1\otimes a$ from $X$ and $A$ respectively into $X\infty A$ are
bialgebra morphisms. Also, $x\infty a=(x\infty 1)(1\infty a)$ for $a\in
A$ and $x\in X$.

\sect{Biperfect Weak Hopf Algebras And Quasi-Bicrossed Products}\label{sec2}

\begin{definition}\label{def2.1}
A weak Hopf algebra $H$ is called: (i) a {\em perfect weak Hopf algebra}
\cite{L2} if its weak antipode $T$ is an anti-bialgebra morphism  satisfying
 $(id\ast T)(H)\subseteq C(H)$ (the center of $H$);
 (ii) a {\em coperfect weak Hopf algebra}
  if its weak antipode $T$ is an anti-bialgebra morphism  satisfying
 $\sum_{(x)}x'T(x'')\otimes x'''=\sum_{(x)}x''T(x''')\otimes x'$ for any
$x\in H$; (iii) a {\em biperfect weak Hopf algebra} if it is perfect and
also coperfect.
\end{definition}

{}From Proposition 1.2 in \cite{L1}, we know that if the weak antipode $T$ of
a weak Hopf algebra $H=(H,m,u,\Delta,\varepsilon,T)$ is an  invertible
anti-algebra morphism, then $H^{op}=(H,m^{op},
u,\Delta,\varepsilon)$ and $H^{cop}=(H,m,u,\Delta^{op},\varepsilon)$ are
both weak Hopf algebras with weak antipode $T^{-1}$.

\begin{lemma}\label{lemma2.1}
Suppose that $H=(H,m,u,\Delta,\varepsilon,T)$ is a weak Hopf algebra with
$T$ invertible, then $H$ is perfect (resp. coperfect) if and only if
$H^{op}$ (resp. $H^{cop}$) is also perfect (resp. coperfect).
\end{lemma}
{\em Proof}:  When $H$ is coperfect, then
$\sum_{(x)}x'T(x'')\otimes x'''=\sum_{(x)}x''T(x''')\otimes x'$ for any
$x\in H$. Thus, \[
(T\otimes 1)\sum_{(x)}(x''T^{-1}(x')\otimes x''')=
(T\otimes 1)\sum_{(x)}(x'''T^{-1}(x'')\otimes x').
\]
It follows that
$\sum_{(x)}(x''T^{-1}(x')\otimes x''')=\sum_{(x)}(x'''T^{-1}(x'')\otimes x')$
since $T$ is invertible. This means that $H^{cop}$ is coperfect on the weak
antipode $T^{-1}$.
It is similar to prove the result in the case that $H$ is perfect. \#

For a finite dimensional weak Hopf algebra $H=(H,m,u,\Delta,\varepsilon,T)$
, we know \cite{L1} that $H^{\ast}=(H^{\ast},\Delta^{\ast},\varepsilon^{\ast},
m^{\ast},u^{\ast},T^{\ast})$ is a weak Hopf algebra with weak antipode
$T^{\ast}$.

\begin{lemma}\label{lemma2.2}
A finite dimensional weak Hopf algebra $H$ is perfect (resp. coperfect)
 if and only if its duality $H^{\ast}$ is coperfect (resp. perfect).
\end{lemma}
{\em Proof}: ``only if": When $H$ is perfect, we need to prove that for
$f\in H^{\ast}$,
\[
\sum_{(f)}f'T^{\ast}(f'')\otimes f'''=\sum_{(f)}f''T^{\ast}(f''')\otimes f'.
\]
In fact, for $a$, $b\in H$,
\begin{eqnarray*}
(\sum_{(f)}f'T^{\ast}(f'')\otimes f''')(a\otimes b)
&=&\sum_{(f)}(f'T^{\ast}(f''))(a)f'''(b)
=\sum_{(f)(a)}f'(a')T^{\ast}(f'')(a'')f'''(b)\\
&=&\sum_{(f)(a)}f'(a')f''(T(a''))f'''(b)=\sum_{(a)}f(a'T(a'')b)\\
&=&\sum_{(a)}f(ba'T(a''))=\sum_{(f)(a)}f'(b)f''(a')f'''(T(a''))\\
&=&\sum_{(f)(a)}f'(b)f''(a')T^{\ast}(f''')(a'')
=\sum_{(f)}f'(b)(f''T^{\ast}(f'''))(a)\\
&=&\sum_{(f)}(f''T^{\ast}(f''')\otimes f')(a\otimes b).
\end{eqnarray*}

When $H$ is coperfect, we need to prove that $(id_{H^{\ast}}\ast T^{\ast})
(H^{\ast})\subseteq C(H^{\ast})$.

In fact, for $f$, $g\in H^{\ast}$, $x\in H$,
\begin{eqnarray*}
(g(id_{H^{\ast}}\ast T^{\ast})(f))(x)
&=&\sum_{(x)}(g(x')(id_{H^{\ast}}\ast T^{\ast})(f))(x'')
=\sum_{(x)}g(x')\Delta^{\ast}(id_{H^{\ast}}\otimes T^{\ast})m^{\ast}(f)(x'')
\nonumber\\
&=&\sum_{(x)}g(x')(id_{H^{\ast}}\otimes T^{\ast})m^{\ast}(f)(x''\otimes x''')
=\sum_{(x)}g(x')m^{\ast}(f)(x''\otimes T(x'''))\nonumber\\
&=&\sum_{(x)}g(x')f(x''T(x'''))
=\sum_{(x)}f(x'T(x''))g(x''')\nonumber\\
&=&((id_{H^{\ast}}\ast T^{\ast})(f)g)(x).
\end{eqnarray*}

It is easy to see that $T^{\ast}$ is an anti-bialgebra morphism from
the same fact of $T$.

``if": It follows from $H\cong H^{\ast\ast}$.
\#

\begin{corollary}\label{cor2.3}
A finite dimensional weak Hopf algebra $H$ is biperfect
 if and only if its duality $H^{\ast}$ is biperfect.
\end{corollary}

\begin{lemma}\label{lemma2.4}
Suppose that $H=(H,m,u,\Delta,\varepsilon,T)$ is a finite dimensional weak
 Hopf algebra and its weak antipode $T$ is an invertible 
 anti-bialgebra morphism. Then
(i) $H$ is perfect if and only if $(T\ast id)(H)\subseteq C(H)$; (ii) $H$
is coperfect if and only if $\sum_{(x)}T(x')x''\otimes x'''
=\sum_{(x)}T(x'')x'''\otimes x'$ for any $x\in H$.
\end{lemma}
{\em Proof}: (i) follows from Lemma 1.1 in \cite{L1}.

(ii) Similar to the proof of Lemma III.2, we can prove that $H$
satisfies $(T\ast id)(H)\subseteq C(H)$ if and only if $H^{\ast}$
satisfies $\sum_{(f)}T^{\ast}(f')f''\otimes
f'''=\sum_{(f)}T^{\ast}(f'')f'''\otimes f'$
 for $f\in H^{\ast}$.

Then $H$ is coperfect if and only if $H^{\ast}$ is perfect, if and only if
 $(T^{\ast}\ast id)(H^{\ast})\subseteq C(H^{\ast})$. But $H\cong(H^{\ast})
^{\ast}$. So, if and only if $\sum_{(x)}T(x')x''\otimes x'''
=\sum_{(x)}T(x'')x'''\otimes x'$ for $x\in H$. \#

The concept of a Hopf pair of Hopf algebras was introduced by
M.Takeuchi in \cite{T}, which plays a valid role in the study of
the theory of quantum groups. Now, we generalize this and introduce some
similar concepts corresponding to the weak Hopf algebras.

\begin{definition}\label{def2.2}
 (i) Suppose that $A$ and $X$ are weak Hopf algebras with weak
antipodes $S_{A}$ and $S_{X}$, respectively.
We call $(X,A)$ a weak Hopf pair, if there
exists a non-singular bilinear form $< , >$ from $X\otimes A$ to $k$
satisfying
\begin{equation}
<x,ab>=\sum_{(x)}<x',a><x'',b>,\label{e9}
\end{equation}
\begin{equation}
<x,1_{A}>=\varepsilon (x),\label{e10}
\end{equation}
\begin{equation}
<xy,a>=\sum_{(a)}<x,a'><y,a''>,\label{e11}
\end{equation}
\begin{equation}
<1_{X},a>=\varepsilon (a),\label{e12}
\end{equation}
\begin{equation}
<S_{X}(x),a>=<x,S_{A}(a)>,\label{e13}
\end{equation}
where $x$, $y\in X$, $a$, $b\in A$.

(ii) In (i), moreover, if  $S_{A}$ is invertible and (\ref{e9}) and (\ref{e13})
are replaced with the following (\ref{e14}) and (\ref{e15}):
\begin{equation}
<x,ab>=\sum_{(x)}< x'',a>< x',b>,\label{e14}
\end{equation}
\begin{equation}
<S_{X}(x),a>=<x,S_{A}^{-1}(a)>,\label{e15}
\end{equation}
respectively, we call  $(X,A)$ a {\em weak Hopf skew-pair}.
\end{definition}

{}From \cite{L1}, $A^{op}=(A,\mu^{op},\eta ,\Delta,\varepsilon
,S_{A}^{-1})$ is a weak Hopf algebra when $S_{A}$ is invertible. Therefore
$(X,A)$ is a weak Hopf skew-pair if and only if $(X,A^{op})$ is a weak Hopf
pair in the case where $S_{A}$ is invertible.

We know from \cite{L4} that for two perfect weak Hopf algebras
$A$ and $X$ with weak antipodes
$S_{A}$ and $S_{X}$ respectively, suppose that $A$ is cocommutative,
$S_{A}$ is invertible and $(X,A)$ is a weak Hopf skew-pair, then $(X,A)$ is
a quasi-matched pair of bialgebra. We want to generalize this result to the
case that $A$ is non-cocommutative. In fact, we have the following lemma:

\begin{lemma}\label{lemma2.5}
For two perfect weak Hopf algebras $A$ and $X$ with weak antipodes
$S_{A}$ and $S_{X}$ respectively, suppose that $S_{A}$ is invertible and
$(X,A)$ is a weak Hopf skew-pair. Then $A$ and $X$ are both biperfect.
\end{lemma}
{\em Proof}: For $x\in X$, $a$, $b\in A$, since $A$ is perfect, we have
\begin{eqnarray*}
\sum_{(x)}<x'S_{X}(x''),a><x''',b>&=&
\sum_{(x)(a)}<x',a'><S(x''),a''><x''',b>\nonumber\\
&=&\sum_{(x)(a)}<x',a'><x'',S_{A}^{-1}(a'')><x''',b>\nonumber\\
&=&\sum_{(a)}<x,bS_{A}^{-1}(a'')a'>
=\sum_{(a)}<x,S_{A}^{-1}(a'')a'b>\nonumber\\
&=&\sum_{(x)(a)}<x',b><x'',a'><x''',S_{A}^{-1}(a'')>\nonumber\\
&=&\sum_{(x)(a)}<x',b><x'',a'><S_{X}(x'''),a''>\nonumber\\
&=&\sum_{(x)}<x''S_{X}(x'''),a><x',b>,
\end{eqnarray*}
and hence $\sum_{(x)}(x'S_{X}(x'')\otimes x''')=
\sum_{(x)}(x''S_{X}(x''')\otimes x')$. It means that $X$ is coperfect.

Similarly, for $x$, $y\in X$, $a\in A$, since $X$ is perfect, we can prove
\[
\sum_{(a)}<x,a''S_{A}^{-1}(a')><y,a'''>=\sum_{(a)}<x,a'''S_{A}^{-1}(a'')><y,a'>.
\]
Hence $\sum_{(a)}(a''S_{A}^{-1}(a')\otimes a''')=
\sum_{(a)}(a'''S_{A}^{-1}(a'')\otimes a')$. Thus $\sum_{(a)}(a'S_{A}(a'')
\otimes a''')=\sum_{(a)}(a''S_{A}(a''')\otimes a')$. It follows that $A$ is
coperfect. \#

\begin{theorem}\label{th2.6}
For two perfect weak Hopf algebras $A$ and $X$ with weak antipodes
$S_{A}$ and $S_{X}$ respectively, suppose that $S_{A}$ is invertible and
$(X,A)$ is a weak Hopf skew-pair. Then $(X,A)$ is a quasi-matched pair of
bialgebras with
\[
a\triangleright x=\sum_{(x)}<x'S_{X}(x'''),a>x'',
\]
\[
a\triangleleft x=\sum_{(a)}<x,S_{A}^{-1}(a''')a'>a''
\]
so as to get a
quasi-bicrossed product $X\infty A$, denoted as $D(X,A)$.
\end{theorem}
{\em Proof}: Firstly, we can verify easily the following:
\begin{equation}
<a\triangleright x,b>=\sum_{(a)}<x,S_{A}^{-1}(a'')ba'>,\label{e16}
\end{equation}
\begin{equation}
<y,a\triangleleft x>=\sum_{(x)}<x'yS_{X}(x''),a>\label{e17}
\end{equation}
for $a$, $b\in A$, $x$, $y\in X$.

Now we prove that $A$ and $X$  are a right $X$-quasi-module coalgebra
and a left $A$-quasi-module coalgebra with the action $\triangleleft$
and $\triangleright$, respectively. In fact,
for any $a\in A$, $x$, $y$, $z\in X$, we
have
\begin{eqnarray*}
<z,a\triangleleft(xy)>&=&\sum_{(xy)}<(xy)''zS_{X}((xy)'),a> \\
&=&\sum_{(x)(y)}<x''y''zS_{X}(y')S_{X}(x'),a>=<z,(a\triangleleft
x)\triangleleft y>,
\end{eqnarray*}
then  $a\triangleleft (xy)=(a\triangleleft x)\triangleleft y$;
$<z,a\triangleleft 1>=<1zS_{X}(1),a>=<z,a>$,
then $a\triangleleft 1=a$. Thus $A$ is a right $X$-module. On the other
hand,
\begin{eqnarray*}
<y\otimes z,\sum_{(a)(x)}(a'\triangleleft x')\otimes
(a''\triangleleft x'')>&=&\sum_{(a)(x)}<y,a'\triangleleft x'>
<z,a''\triangleleft x''>\\
&=&\sum_{(a)(x)}<x'yS_{X}(x''),a'><x'''zS_{X}(x^{(4)}),a''>\\
&=&\sum_{(x)}<x'yS_{X}(x'')x'''zS_{X}(x^{(4)}),a>\\
&=&\sum_{(x)}<x'yzS_{X}(x''),a>\\
&=&<yz,a\triangleleft x>=<y\otimes z,\Delta(a\triangleleft x)>,
\end{eqnarray*}
then $\Delta (a\triangleleft x)=\sum_{(a)(x)}(a'\triangleleft
x')\otimes (a''\triangleleft x'')$. It means that $A$ is a right
$X$-quasi-module-coalgebra.  Similarly, we get that $X$ is a left
$A$-quasi-module-coalgebra.

Moreover, we can see that (\ref{e2}) and (\ref{e4}) are trivial according
to (\ref{e12}) and
(\ref{e10}) and the definition of $\triangleright$  and $\triangleleft$ . And,
 using of Lemma \ref{lemma2.5}, we have
\begin{eqnarray*}
<\sum_{(a)(x)}(a'\triangleright x')((a''\triangleleft
x'')\triangleright y),b>
&=&<\sum_{(a)(x)}(a'\triangleright
x')(<x'',S_{A}^{-1}(a^{(4)})a''>a'''\triangleright y),b>\\
&=&\sum_{(a)(b)(x)}<x'',S_{A}^{-1}(a^{(4)})a''><(a'\triangleright
x'),b'><(a'''\triangleright y),b''>\\
&=&\sum_{(a)(b)(x)}<x'',S_{A}^{-1}(a^{(5)})a'''>\\
& &~~~~~~<x',S_{A}^{-1}(a'')b'a'> <(a^{(4)}\triangleright y),b''>\\
&=&\sum_{(a)(b)}<x,S_{A}^{-1}(a^{(5)})a'''S_{A}^{-1}(a'')b'a'>
<(a^{(4)}\triangleright y),b''>\\
&=&\sum_{(a)(b)}<x,S_{(A)}^{-1}(a''')b'a'><(a''\triangleright y),b''>\\
&=&\sum_{(a)(b)}<x,S_{A}^{-1}(a^{(4)})b'a'><y,S_{A}^{-1}(a''')b''a''>\\
&=&\sum_{(a)}<xy,S_{A}^{-1}(a'')ba'>=<a\triangleright (xy),b>,
\end{eqnarray*}
then $a\triangleright (xy)=\sum_{(a)(x)}(a'\triangleright
x')((a''\triangleleft x'')\triangleright y)$, i.e. eq.(\ref{e1}) holds.
Similarly, we get that $(ab)\triangleleft x=\sum_{(b)(x)}(a\triangleleft
(b'\triangleright x'))(b''\triangleleft x'')$, i.e. eq.(\ref{e3}) holds.
Moreover,
\begin{eqnarray*}
\sum_{(a)(x)}(a'\triangleleft x')\otimes (a''\triangleright x'')
&=&\sum_{(a)(x)}<x',S_{A}^{-1}(a''')a'>a''\otimes
<x''S_{X}(x^{(4)}),a^{(4)}>x'''\\
&=&\sum_{(a)(x)}<x',S_{A}^{-1}(a''')a'><x'',a^{(4)}><
S_{X}(x^{(4)}),a^{(5)}>a''\otimes x'''\\
&=&\sum_{(a)(x)}<x',a^{(4)}S_{A}^{-1}(a''')a'><S_{X}(x'''),a^{(5)}>a''\otimes
x''\\
&=&\sum_{(a)(x)}<x',a'''S_{A}^{-1}(a'')a'><S_{X}(x'''),a^{(5)}>a^{(4)}\otimes
x''\\
&=&\sum_{(a)(x)}<x',a'>< S_{X}(x'''),a'''>a''\otimes
x''\\
&=&\sum_{(a)(x)}<x',a'>< x''',S_{A}^{-1}(a''')>a''\otimes x''\\
&=&\sum_{(a)(x)}<x',a'>< x''',S_{A}^{-1}(a^{(5)})a^{(4)}S_{A}^{-1}(a''')
>a''\otimes x''\\
&=&\sum_{(a)(x)}<x',a'>< x''',S_{A}^{-1}(a^{(5)})a'''
S_{A}^{-1}(a'')>a^{(4)}\otimes x''\\
&=&\sum_{(a)(x)}<x',a'><x''',S_{A}^{-1}(a'')><x^{(4)},
S_{A}^{-1}(a^{(5)})a'''>a^{(4)}\otimes x''\\
&=&\sum_{(a)(x)}<x',a'><S_{X}(x'''),a''>
<x^{(4)},S_{A}^{-1}(a^{(5)})a'''>a^{(4)}\otimes x''\\
&=&\sum_{(a)(x)}<x'S_{X}(x'''),a'><x^{(4)},S_{A}^{-1}(a^{(4)})a''>a'''
\otimes x''\\
&=&\sum_{(a)(x)}(a''\triangleleft x'')\otimes (a'\triangleright x'),
\end{eqnarray*}
then eq.(\ref{e5}) holds. In a word, $(X,A)$ is a quasi-matched pair of
bialgebras. Hence we get a quasi-bicrossed product $X\infty A$, denoted also
as $D(X,A)$.   \#

Note that in this theorem, $A$ is not required to be cocommutative. So, this
theorem is a big improvement on the result obtained in \cite{L4}.

\sect{Quantum Doubles Of Biperfect Weak Hopf Algebras}

In Theorem \ref{th2.6}, when $A=H$ is a finite dimensional biperfect weak
Hopf algebra
with invertible weak antipode $T$, we set $X=H^{\ast cop}$ and suppose that
 $< , >$ is the bilinear form of $H$ and its dual
$H^{\ast}$ as linear spaces. It was known in \cite{L1} that $(H^{\ast cop}=
(H^{\ast},\Delta^{\ast},\varepsilon^{\ast},
(m^{\ast})^{op},u^{\ast},(T^{\ast -1})$. It is easy to see that
$(H^{\ast cop},H)$ is a weak Hopf skew-pair. Then $(H^{\ast cop},H)$ is a
 quasi-matched pair of bialgebras with $a\triangleright f=
\sum_{(f)}<f'T^{\ast -1}(f'''),a>f''$ and $a\triangleleft
f=\sum_{(a)}<f,T^{-1}(a''')a'>a''$  for $a\in H$ and $f\in H^{\ast
cop}$ so as to get a quasi-bicrossed product $D(H^{\ast cop},H)
=H^{\ast cop}\infty H$, denoted briefly as $D(H)$ and called the
{\em quantum
 double} of $H$.

\begin{proposition}\label{prop3.1}
Let $H=(H,m,u,\Delta,\varepsilon,T)$ be a finite dimensional biperfect weak
Hopf algebra with invertible $T$. Then the multiplication in $D(H)
=H^{\ast cop}\infty H$ is given by
\[
(f\infty a)(g\infty b)=\sum_{(a)}fg(T^{-1}(a''')?a')\infty a''b
\]
for $f$, $g\in H^{\ast cop}$, $a$, $b\in H$.
\end{proposition}
{\em Proof}:
\begin{eqnarray*}
(f\infty a)(g\infty b)&=&\sum_{(a)(g)}f(a'\triangleright g')\infty
(a''\triangleleft g'')b\\
&=&\sum_{(a)(g)}fg'(T^{-1}(a'')?a')\infty g''(T^{-1}(a^{(5)})a''')a^{(4)}b\\
&=&\sum_{(a)}fg(T^{-1}(a^{(5)})a'''T^{-1}(a'')?a')\infty a^{(4)}b\\
&=&\sum_{(a)}fg(T^{-1}(a^{(5)})?a'''T^{-1}(a'')a')\infty a^{(4)}b
=\sum_{(a)}fg(T^{-1}(a''')?a')\infty a''b.~~~ \#
\end{eqnarray*}

Now we have the following main result:
\begin{theorem}\label{th3.2}
Let $H=(H,m,u,\Delta,\varepsilon,T)$ be a finite dimensional biperfect
weak Hopf algebra with invertible $T$.
Then the quantum double $D(H)$ of $H$ is quasi-braided equipped with a
 quasi-R-matrix $R=\sum_{i\in I}(1\infty e_{i})\otimes(e^{i}\infty 1)\in
 D(H)\otimes D(H)$ where $\{e_{i}\}_{i\in I}$ is a basis of the $k$-vector
 space $H$ together with its dual basis $\{e^{i}\}_{i\in I}$ in $H^{\ast cop}$.
Hence $R$ is a solution  of the quantum Yang-Baxter equation.
\end{theorem}
{\em Proof}: For $f\in H^{\ast cop}$, $a\in H$,
\begin{eqnarray*}
\Delta^{op}(f\infty a)R&=&\sum_{i\in I}\sum_{(f)(a)}(f''\infty a'')(1\infty
e_{i})\otimes(f'\infty a')(e^{i}\infty 1)\\
&=&\sum_{i\in I}\sum_{(f)(a)}(f''\varepsilon(T^{-1}(a^{(6)})?a^{(4)})\infty
a^{(5)}e_{i})\otimes
(f'e^{i}(T^{-1}(a''')?a')\infty a'')\\
&=&\sum_{i\in I}\sum_{(f)(a)}(f''\varepsilon(T^{-1}(a^{(6)}))\varepsilon
(a^{(4)})\varepsilon\infty  a^{(5)}e_{i})\otimes (f'e^{i}(T^{-1}(a''')?a')
\infty a'')\\
&=&\sum_{i\in I}\sum_{(f)(a)}(f''\infty  a^{(4)}e_{i})\otimes (f'e^{i}(T^{-1}
(a''')?a')\infty a'');
\end{eqnarray*}
\begin{eqnarray*}
R\Delta(f\infty a)&=&\sum_{i\in I}\sum_{(f)(a)}(\varepsilon\infty e_{i})(f'
\infty a')\otimes(e^{i}\infty 1)(f''\infty a'')\\
&=&\sum_{i\in I}\sum_{(f)(a)(e_{i})}(\varepsilon f'(T^{-1}(e_{i}''')?e_{i}')
\infty e_{i}''a')\otimes(e^{i}f''(T^{-1}(1)?1)\infty 1a'')\\
&=&\sum_{i\in I}\sum_{(f)(a)(e_{i})}(f'(T^{-1}(e_{i}''')?e_{i}')\infty
 e_{i}''a')\otimes(e^{i}f''\infty a'').
\end{eqnarray*}

For every $b$, $c\in H$, $u$, $v\in H^{op\ast}$, let $\xi=b\otimes u\otimes
 c\otimes v$. Then
\begin{eqnarray*}
<\Delta^{op}(f\infty a)R,\xi>&=&\sum_{i\in I}
\sum_{(f)(a)}f''(b)u(a^{(4)}e_{i})(f'e^{i}(T^{-1}(a''')?a'))(c)v(a'')\\
&=&\sum_{i\in I}\sum_{(f)(a)(c)}f''(b)u(a^{(4)}e_{i})f'(c')e^{i}(T^{-1}(a''')c''a')v(a'')\\
&=&\sum_{i\in I}\sum_{(f)(a)(c)}f''(b)u(a^{(4)}e^{i}
(T^{-1}(a''')c''a')e_{i})f'(c')v(a'')\\
&=&\sum_{(a)(c)}f(bc')u(a^{(4)}T^{-1}(a''')c''a')v(a'')\\
&=&\sum_{(a)(c)}f(bc')u(c''a^{(4)}T^{-1}(a''')a')v(a'')\\
&=&\sum_{(a)(c)}f(bc')u(c''a'''T^{-1}(a'')a')v(a^{(4)})
=\sum_{(a)(c)}f(bc')u(c''a')v(a'');
\end{eqnarray*}
\begin{eqnarray*}
<R\Delta(f\infty a),\xi>&=&
\sum_{i\in I}\sum_{(f)(a)(e_{i})}f'(T^{-1}(e_{i}''')be_{i}')u(e_{i}''a')
(e^{i}f'')(c)v(a'')\\
&=&\sum_{i\in I}\sum_{(f)(a)(e_{i})(c)}f'(T^{-1}(e_{i}''')be_{i}')u(e_{i}''a')
e^{i}(c')f''(c'')v(a'')\\
&=&\sum_{i\in I}\sum_{(a)(e_{i})(c)}f(c''T^{-1}(e_{i}''')be_{i}')u(e_{i}''a')
e^{i}(c')v(a'')\\
&=&\sum_{(a)(c)}f(c^{(4)}T^{-1}(c''')bc')u(c''a')v(a'')\\
& &({\rm since}~ \sum_{i\in I}\sum_{(e_{i})}e^{i}(c)e_{i}'
\otimes e_{i}''\otimes e_{i}'''=\sum_{(c)}c'\otimes c''\otimes c''')\\
&=&\sum_{(a)(c)}f(c'''T^{-1}(c'')bc')u(c^{(4)}a')v(a'')\\
&=&\sum_{(a)(c)}f(bc'''T^{-1}(c'')c')u(c^{(4)}a')v(a'')\\
&=&\sum_{(a)(c)}f(bc')u(c''a')v(a'')=<\Delta^{op}(f\infty a)R,\xi>.
\end{eqnarray*}

Therefore $\Delta^{op}(f\infty a)R=R\Delta(f\infty a)$. Then $H^{\ast cop}
\infty H$ is an almost quasi-cocommutative almost bialgebra with a universal
quasi-R-matrix $R$.

We can prove in a way similar to the proof of Theorem 2.11 in \cite{L1} that
\[
(\Delta\otimes id_{H})(R)=R_{13}R_{23}; (id_{H}\otimes\Delta)(R)=R_{13}R_{12}.
\]

It means that $H^{op\ast}\infty H$ is quasi-braided. Thus, by Proposition
2.8 in \cite{L1}, $R$ is a quasi-R-matrix.
\#

Note that since a cocommutative weak Hopf algebra must be
coperfect, it means that Theorem IV.2 is a generalization
of the one  in \cite{L1} on the quantum double of a finite dimensional
cocommutative perfect weak Hopf algebra.

 It is easy to see that for a finite Clifford monoid
$S=\{s_{1},...,s_{n}\}$ (see \cite{P} and \cite{L1}),
 $kS$ is a finite dimensional biperfect weak Hopf algebra with invertible
weak antipode $T_{S}$ satisfying $T_{S}(s)=s^{-1}$ for $s\in S$ . Then by
Theorem \ref{th3.2}, the quantum double $D(kS)$ is quasi-braided equipped with a
quasi-R-matrix $R=\sum_{i=1}^{n}(1\infty s_{i})\otimes(s_{i}^{\ast}\infty 1)
\in D(kS)\otimes D(kS)$ where $s_{i}^{\ast}$ is the duality of $s_{i}$ in
$(kS)^{\ast cop}$. Thus, $R$ is a solution  of the quantum Yang-Baxter
equation. But, this is also an example of a quantum quasi-double in Theorem
 2.11 of \cite{L1} constructed
 from a finite dimensional cocommutative perfect weak Hopf
algebra. So, it is very necessary to find an example of the quantum double
from a finite dimensional biperfect
 weak Hopf algebra with invertible weak antipode which is not cocommutative.

For two biperfect weak Hopf algebras $H$ and $K$, it is easy to
prove that the tensor product $H\otimes K$ is also a biperfect weak Hopf
algebra with the comultiplication $\Delta=(I\otimes T\otimes I)(\Delta_{H}
\otimes\Delta_{K})$, the multiplication $m=(m_{H}\otimes m_{K})(I\otimes T
\otimes I)$, the unit $1=1_{H}\otimes 1_{K}$, the counit $\varepsilon=
\varepsilon_{H}\otimes\varepsilon_{K}$ and the weak antipode $T=T_{H}\otimes
T_{K}$. $H\otimes K$ is commutative (resp. cocommutative) if and only if $H$
 and $K$ are so.

For a finite non-commutative Clifford monoid $S$, let $H=kS$ with the weak
antipode $T_{S}$, then
 $K=(kS)^{\ast}$ is also a finite dimensional biperfect weak Hopf algebra
with invertible weak antipode $T_{S}^{\ast}$. Thus, we get a finite
 dimensional biperfect weak Hopf algebra $A=kS\otimes(kS)^{\ast}$ with
 invertible weak
 antipode $T=T_{S}\otimes T_{S}^{\ast}$, which is indeed not a Hopf algebra
unless $S$ is a group. Since $kS$ is non-commutative,
 $(kS)^{\ast}$ is non-cocommutative. Hence $A$ is non-commutative and
non-cocommutative.  By Theorem \ref{th3.2}, the quantum double $D(A)$ of $A$ is
quasi-braided equipped with a quasi-R-matrix as a solution of the quantum
Yang-Baxter equation.  This construction is different from
that of Theorem 2.11 in \cite {L1}. It implies that in Theorem \ref{th3.2}, the
quantum double of a finite dimensional biperfect weak Hopf algebra is indeed
a generalization of that of a finite dimensional Hopf algebra.

We know in \cite{K} that the R-matrix of quantum double of a finite
dimensional Hopf algebra is invertible. But, for the quasi-R-matrix in
Theorem \ref{th3.2}, we can only get the regularity as following:

\begin{proposition}\label{prop3.3}
For a finite dimensional biperfect weak Hopf algebra $H$ with invertible
weak antipode $T$, the quasi-R-matrix
$R=\sum_{i=1}^{n}(\varepsilon\infty e_{i})\otimes(e_{i}^{\ast}\infty
1)$ of its quantum quasi-double $D(H)$ is a von Neumann regular element in
$D(H)\otimes D(H)$ with its inverse $\bar{R}=\sum_{i=1}^{n}(\varepsilon\infty
e_{i})\otimes(e_{i}^{\ast}T\infty 1)$ where $\{e_{1}, ... , e_{n}\}$ is a
basis of $H$ and $\{e_{1}^{\ast}, ... , e_{n}^{\ast}\}$ is the dual
basis in $H^{\ast}$.
\end{proposition}
{\em Proof}:  For any $\xi=b\otimes u\otimes c\otimes v\in H\otimes
H^{\ast}\otimes H\otimes H^{\ast}$, we have
\begin{eqnarray*}
<R\bar{R}R,\xi>&=&\sum_{i,j,l}\varepsilon(b)u(e_{i}e_{j}e_{l})
(e_{i}^{\ast}(e_{j}^{\ast}T)e_{l}^{\ast})(c)v(1) \\
&=&\varepsilon(b)v(1)\sum_{(c)}u(\sum_{i=1}^{n}e_{i}e_{i}^{\ast}(c')
\sum_{j=1}^{n}e_{j}e_{j}^{\ast}(T(c''))\sum_{l=1}^{n}e_{l}e_{l}^{\ast}(c'''))
\\
&=&\varepsilon(b)v(1)u(\sum_{(c)}c'T(c'')c''')
=\varepsilon(b)v(1)u(c) \\
&=&\varepsilon(b)v(1)\sum_{i=1}^{n}u(e_{i})e_{i}^{\ast}(c)
=<R,\xi>,
\end{eqnarray*}
then $R\bar{R}R=R$. Similarly, $<\bar{R}R\bar{R},\xi>=<\bar{R},\xi>$,
then $\bar{R}R\bar{R}=\bar{R}$.   \#

For a left $D(H)$-module $V$, define $C_{V,V}^{R}$ satisfying
$C_{V,V}^{R}(v\otimes w)=\tau(R(v\otimes w))$
for $v$, $w\in V$, then $C_{V,V}^{R}$ is a solution of the classical
Yang-Baxter equation (see \cite{L2}) where $\tau$ is the flip
map defined as $\tau(v_{1}\otimes v_{2})=v_{2}\otimes v_{1}$. From
Propositon \ref{prop3.3}, it is easy to prove that:

\begin{corollary}\label{cor3.4}
For a finite dimensional biperfect weak Hopf algebra $H$
with invertible weak antipode $T$, let $V$ be a left $D(H)$-module. Then
$C_{V,V}^{R}$ is regular in the endomorphism monoid of $V\otimes V$,
with its inverse $C_{V,V}^{\bar{R}}$ satisfying
$C_{V,V}^{\bar{R}}(v\otimes w)=\tau(\bar{R}(v\otimes w))$ for $v$, $w\in V$.
\end{corollary}

Now we discuss the representation-theoretic interpretation of $D(H)$.

\begin{definition}\label{def3.1}
(see \cite{K})  For a bialgebra $H$ over $k$, a {\em crossed
$H$-bimodule} $V$ is a vector space together with linear maps ${\mu}_{V}:
H\otimes V\longrightarrow V$  and ${\Delta}_{V}: V\longrightarrow V
\otimes H$ such that

(i) the map ${\mu}_{V}$ and ${\Delta}_{V}$ turn $V$ into a left
$H$-module and a right $H$-comodule respectively;

(ii) $\sum_{(a)(\beta)}a'\beta_{V}\otimes
a''\beta_{H}=\sum_{(a)(a''\beta)}(a''\beta)_{V}\otimes (a''\beta)_{H}a'$
 for all $a\in H$ and $\beta\in V$ where we set ${\mu}_{V}(a\otimes
\beta)=a\beta$ and
$\Delta_{V}(\beta)=\sum_{(\beta)}\beta_{V}\otimes\beta_{H}$.
\end{definition}

\begin{theorem}\label{th3.5}
Suppose $H$ is a finite dimensional biperfect weak Hopf algebra
with invertible weak antipode $T$. Then for a $k$-linear space
$V$, the following statements are equivalent:

(i) $V$ is a left $D(H)$-module;

(ii) $V$ is a crossed $H$-bimodule $V$ and satisfies
\begin{equation}
\sum_{(a)(\beta)}T^{-1}(a''')a''\beta_{H}\otimes
a'\beta_{V}=\sum_{(\beta)}\beta_{H}\otimes a\beta_{V}\label{e18}
\end{equation}
for all $a\in H$ and $\beta\in V$.
\end{theorem}
{\em Proof}: ``$(i)\Longrightarrow (ii)$":  Let $V$ be a left
$D(H)$-module. Since $H^{\ast cop}\cong H^{\ast cop}\infty 1$ and
$H\cong\varepsilon\infty H$ are subalgebras of $D(H)$, $V$ is a
left $H$-module and also a left $H^{\ast}$-module with
$a\beta=(1\infty a)\beta$ and $x\beta=(x\infty 1)\beta$  for $a\in
H$, $x\in H^{\ast}$, $\beta\in V$. Then, $(ax)(\beta)=a(x\beta)$.
But, by Proposition \ref{prop3.1}, we get
\begin{eqnarray*}
ax&=&(1\infty a)(x\infty 1)=\sum_{(a)(x)}<x,T^{-1}(a''')?a'>\infty a''\\
&=&\sum_{(a)(x)}(<x,T^{-1}(a''')?a'>\infty 1)(1\infty a'').
\end{eqnarray*}
Hence, $a(x\beta)=\sum_{(a)(x)}<x,T^{-1}(a''')?a'>(a''\beta)$.

One must show that $V$ can be endowed with a crossed $H$-bimodule
structure. For $\mu_{V}$, we define $\mu_{V}(a\otimes\beta)=a\beta$.

Given a basis $\{e_{1}, ... , e_{n}\}$ of $H$ and the dual basis
$\{e_{1}^{\ast}, ... , e_{n}^{\ast}\}$ of
$H^{\ast}$, note that $x=\sum_{i=1}^{n}x(e_{i})e_{i}^{\ast}$ and
$a=\sum_{i=1}^{n}e_{i}^{\ast}(a)e_{i}$ for $x\in H^{\ast}$, $a\in H$.

 Define $\Delta_{V}: V\longrightarrow V\otimes H$ satisfying
$\Delta_{V}(\beta)=\sum_{i}e_{i}^{\ast}\beta\otimes e_{i}$ for any
$\beta\in V$. Consider the dual $\Delta_{V}^{\ast}$ of
$\Delta_{V}$, we have that for any $\alpha\in V^{\ast}$, $\beta\in
V$, $x\in H^{\ast}$,
\begin{eqnarray*}
<\Delta_{V}^{\ast}(\alpha\otimes x),\beta>&=&<\alpha\otimes
x,\Delta_{V}(\beta)>=\sum_{i=1}^{n}<\alpha,e_{i}^{\ast}\beta><x,e_{i}>\\
&=&<\alpha,(\sum_{i=1}^{n}<x,e_{i}>e_{i}^{\ast})\beta>=<\alpha,x\beta>;
\end{eqnarray*}
in particular, for $x=\varepsilon$  (the identity of
$H^{\ast}$) and any $\beta\in V$,
\[
<\Delta_{V}^{\ast}(\alpha\otimes\varepsilon),\beta>=<\alpha,\beta>.
\]
Then $\Delta_{V}^{\ast}(\alpha\otimes\varepsilon)=\alpha$. It follows that,
for any $y\in H^{\ast}$,
\begin{eqnarray*}
<\Delta_{V}^{\ast}(\Delta_{V}^{\ast}(\alpha\otimes x)\otimes
y),\beta>&=&<\Delta_{V}^{\ast}(\alpha\otimes x),y\beta>\\
&=&<\alpha,x(y\beta)>=<\Delta_{V}^{\ast}(\alpha\otimes xy),\beta>,
\end{eqnarray*}
then $\Delta_{V}^{\ast}(\Delta_{V}^{\ast}(\alpha\otimes x)\otimes
y)=\Delta_{V}^{\ast}(\alpha\otimes xy)$. Hence $V^{\ast}$ is a right
$H^{\ast}$-module under the action $\Delta_{V}^{\ast}$. Therefore, dually,
$V$ becomes a right $H$-comodule under the coaction $\Delta_{V}$.

For $a\in H$, $\beta\in V$, $x\in H^{\ast cop}$, we have
\begin{eqnarray*}
(id\otimes x)(\sum_{(a)(\beta)}a'\beta_{V}\otimes
a''\beta_{H})&=&\sum_{(a)}\sum_{i=1}^{n}a'(e_{i}^{*}\beta)x(a''e_{i})\\
&=&\sum_{(a)(x)}\sum_{i=1}^{n}a'(e_{i}^{*}\beta)x'(e_{i})x''(a'')
=\sum_{(a)(x)}a'(x'\beta)x''(a'')\\
&=&\sum_{(a)(x)}x''(a^{(4)})x'(T^{-1}(a''')?a')(a''\beta)\\
&=&\sum_{(a)}x(a^{(4)}T^{-1}(a''')?a')(a''\beta)\\
&=&\sum_{(a)}x(?a'''T^{-1}(a'')a')(a^{(4)}\beta)=
\sum_{(a)}x(?a')(a''\beta)\\
&=& \sum_{(a)(x)}x'(a')x''(a''\beta)
= \sum_{(a)(x)}\sum_{i=1}^{n}x'(a')x''(e_{i})e_{i}^{*}(a''\beta)\\
&=& \sum_{(a)}\sum_{i=1}^{n}x(e_{i}a')e_{i}^{*}(a''\beta)
=(id\otimes x)(\sum_{(a)}\sum_{i=1}^{n}e_{i}^{*}(a''\beta)\otimes
e_{i}a')\\
&=&(id\otimes x)(\sum_{(a''\beta)(a)}(a''\beta)_{V}\otimes
(a''\beta)_{H}a').
\end{eqnarray*}

It follows that
\begin{equation}
\sum_{(a)(\beta)}a'\beta_{V}\otimes
a''\beta_{H}=\sum_{(a''\beta)(a)}(a''\beta)_{V}\otimes
(a''\beta)_{H}a'.\label{e19}
\end{equation}
 Hence by Definition \ref{def3.1}, $V$ is a crossed $H$-bimodule.

Now, we prove the formula (IV.1). With the $\mu_{V}$ and
$\Delta_{V}$ as defined above, we have proved that for any
$\alpha\in V^{\ast}$, $\beta\in V$, $x\in H^{\ast}$,
$<\Delta_{V}^{\ast}(\alpha\otimes x),\beta>=<\alpha,x\beta>$. But,
 $$<\Delta_{V}^{\ast}(\alpha\otimes x),\beta>=<\alpha\otimes
 x,\beta_{V}\otimes\beta_{H}>=<\alpha,\beta_{V}><x,\beta_{H}>
 =<\alpha,<x,\beta_{H}>\beta_{V}>.$$
 So, it follows
\begin{equation}
x\beta=<x,\beta_{H}>\beta_{V}.\label{e20}
\end{equation}

And, since $V$ is a $D(H)$-module, we have $a(x\beta)=(ax)\beta$.
However, by the formula (IV.3),
$a(x\beta)=\sum_{(\beta)}<x,\beta_{H}>a\beta_{V}$; and by the
formulas (IV.3) and (IV.2),
\begin{eqnarray*}
(ax)\beta&=&(\sum_{(a)}<x,T^{-1}(a''')?a'>\infty
a'')\beta=(\sum_{(a)}<x,T^{-1}(a''')?a'>a'')\beta\\
&=&\sum_{(a)}<x,T^{-1}(a''')?a'>(a''\beta)
=\sum_{(a)(x)}<x''',T^{-1}(a''')><x',a'>x''(a''\beta)\\
&=&\sum_{(a)(x)(a''\beta)}<x''',T^{-1}(a''')><x'',(a''\beta)_{H}>
<x',a'>(a''\beta)_{V}\\
&=&\sum_{(a)(a''\beta)}<x,T^{-1}(a''')(a''\beta)_{H}a'>(a''\beta)_{V}\\
&=&\sum_{(a)(\beta)}<x,T^{-1}(a''')a''\beta_{H}>a'\beta_{V}.
\end{eqnarray*}
Hence, for any $x\in H^{\ast cop}$,
$\sum_{(\beta)}<x,\beta_{H}>a\beta_{V}
=\sum_{(a)(\beta)}<x,T^{-1}(a''')a''\beta_{H}>a'\beta_{V}$.
Then the formula (IV.1) follows.

 ``$(ii)\Longrightarrow (i)$":  Say that $V$ is a crossed
$H$-bimodule about $\mu_{V}$ and $\Delta_{V}$. Then, $V$ is a left
$H$-module about $\mu_{V}$ and a right $H$-comodule about
$\Delta_{V}$. Write $\mu_{V}(a\otimes\beta)=a\beta$ for $a\in H$,
$\beta\in V$. For $x\in H^{\ast}$, $\beta\in V$, let
$x\beta=\sum_{(\beta)}<x,\beta_{H}>\beta_{V}$, where
$\Delta_{V}(\beta)=\sum_{(\beta)}{\beta}_{V}\otimes\beta_{H}$.
Since $\Delta_{V}$ is a right coaction, it is easy to show that
$(xy)\beta=x(y\beta)$ for $y\in H^{\ast}$. Then it follows that
$V$ is a left $H^{\ast cop}$-module.

 Set $(xa)\beta=x(a\beta)$ for $x\in H^{\ast cop}$, $a\in H$,
$\beta\in V$. Then, by (\ref{e18}),
\begin{eqnarray*}
a(x\beta)=\sum_{(\beta)}<x,\beta_{H}>a\beta_{V}=
\sum_{(a)(\beta)}<x,T^{-1}(a''')a''\beta_{H}>a'\beta_{V}=(ax)\beta.
\end{eqnarray*}
where the first equality follows from (VI.3), the second from
(VI.1) and the third from (VI.3) and (VI.2) as proved in
``$(i)\Longrightarrow (ii)$".

Therefore,  $V$ becomes a left $D(H)$-module since $H$ and
$H^{\ast cop}$ are subalgebras of $D(H)$ and the multiplication of
$D(H)$ is determined by the interaction of $H$ and $H^{\ast cop}$.
$\;\;$ \#

\sect{Examples From Matrix Groups}

Now, we give some examples from a concrete Clifford monoid.
The definition of a Clifford semigroup/monoid can be found in \cite{P}
and \cite{L1}.

Let $Y=\{\alpha,\beta,\gamma,\rho,\sigma,\delta\}$ be the
semilattice with multiplication ``$\cdot$" given by the following
table:

\[\begin{tabular}{c|cccccc}

$\cdot$ &$\alpha $& $\beta $ & $\gamma $&$\rho $& $\sigma $&$\delta$\\
\hline
$\alpha $&$\alpha $&$\alpha $&$\alpha $& $\alpha $&$\alpha $&$\alpha$\\
$\beta $&$\alpha $&$\beta$&$\beta $& $\alpha $&$\beta $&$\beta$\\
$\gamma$&$\alpha $&$\beta$&$\gamma $& $\alpha $&$\beta $&$\gamma$\\
$\rho $&$\alpha $&$\alpha $&$\alpha $& $\rho$&$\rho$&$\rho$\\
$\sigma  $&$\alpha $&$\beta$&$\beta $& $\rho $&$\sigma $&$\sigma$\\
$\delta  $&$\alpha $&$\beta$&$\gamma $& $\rho $&$\sigma $&$\delta$\\

\end{tabular}
\]

The partial order in the semilattice $Y$ can be presented as the diagram below:

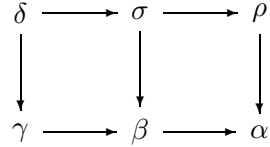
\[
\begin{diagram}
\node{\delta}\arrow{s,r}{}\arrow{e,t}{}\node[1]{\sigma
}\arrow{s,r}{}\arrow{e,t}{}\node[1]{\rho }
\arrow{s,r}{}\\
\node{\gamma
}\arrow{e,t}{}\node[1]{\beta}\arrow{e,t}{}\node[1]{\alpha}
\end{diagram}
\]
Obviously, $\delta$ is the identity of $Y$.

For a ring $R$ with identity, $R^{2\times2}$ denotes the
$2\times2$ full matrix ring over $R$, $U(R)$ the group consisting
of all units in $R$. Let $Z$ be the integer number ring. For a
prime number $p$,  $Z_{p}$ is a field and $U(Z_{p}^{2\times2})$ is
just the $2\times2$ general linear group $GL_{2}(Z_{p})$ over
$Z_{p}$.  Assume that $G_{\alpha}=\{e_{\alpha}\}$ and
$G_{\delta}=\{e_{\delta}\}$ are the trivial groups,
$G_{\beta}=GL_{2}(Z_{2})$, $G_{\gamma}=U(Z_{4}^{2\times2})$,
$G_{\rho}=GL_{2}(Z_{3})$, $G_{\sigma}=U(Z_{6}^{2\times2})$.  Then
$G_{u}\cap G_{v}=\emptyset$ for any $u,v\in Y$, $u\not=v$. Set
$S=\cup_{u\in Y}G_{u}$. We will define a multiplication on $S$
such that $S=\cup_{u\in Y}G_{u}$ becomes a Clifford monoid related
to the semilattice $Y$.

Firstly, we mention the fact that over a commutative ring $R$ with identity,
an $m\times m$ matrix $X$ is invertible if and only if det$X$ is a unit
in $R$.

  Then, for $n=2, 3, 4, 6$, $X=
\left[\begin{array}{cc}
x&y\\
a&b
\end{array}
\right] \in U(Z_{n}^{2\times2})$ if and only if det$X=xb-ay\in
U(Z_{n})$. It is easy to see $U(Z_{6})=\{\bar{1}, \bar{5}\}$,
$U(Z_{4})=\{\bar{1}, \bar{3}\}$, $U(Z_{3})=\{\bar{1}, \bar{2}\}$,
$U(Z_{2})=\{\bar{1}\}$.

A ring homomorphism $\pi_{\sigma,\rho }: Z_{6}\longrightarrow Z_{3}$ can be
defined which satisfies
$\pi_{\sigma,\rho }(\bar{0})=\bar{0}$, $\pi_{\sigma,\rho }(\bar{1})=\bar{1}$,
$\pi_{\sigma,\rho }(\bar{2})=\bar{2}$,
$\pi_{\sigma,\rho }(\bar{3})=\bar{0}$, $\pi_{\sigma,\rho }(\bar{4})=\bar{1}$
and $\pi_{\sigma,\rho }(\bar{5})=\bar{2}$.

For $X=
\left[\begin{array}{cc}
x&y\\
a&b
\end{array}
\right] \in U(Z_{6}^{2\times2})=G_{\sigma}$, we have
det$X=xb-ay=\bar{1}$, or $ \bar{5}$, then $$\pi_{\sigma,\rho
}(x)\pi_{\sigma,\rho }(b)-\pi_{\sigma,\rho }(a) \pi_{\sigma,\rho
}(y)=\bar{1}, \; or\; \bar{2}.$$ It follows
$\left[\begin{array}{cc}
\pi_{\sigma,\rho }(x)&\pi_{\sigma,\rho }(y)\\
\pi_{\sigma,\rho }(a)&\pi_{\sigma,\rho }(b)
\end{array}
\right]
\in GL_{2}(Z_{3})=G_{\rho}$.  Thus, we can expand $\pi_{\sigma,\rho }$ to
make it a group homomorphism from
$G_{\sigma}$ to $G_{\rho }$. For this, it is enough to define
$\pi_{\sigma,\rho }: G_{\sigma}\longrightarrow G_{\rho}$ satisfying
$$\pi_{\sigma,\rho }\left[\begin{array}{cc}
x&y\\
a&b
\end{array}
\right]=
\left[\begin{array}{cc}
\pi_{\sigma,\rho }(x)&\pi_{\sigma,\rho }(y)\\
\pi_{\sigma,\rho }(a)&\pi_{\sigma,\rho }(b)
\end{array}
\right]$$
since  $
\pi_{\sigma,\rho }
\Big(
\left[\begin{array}{cc}
x_{1}&y_{1}\\
a_{1}&b_{1}
\end{array}
\right]\left[\begin{array}{cc}
x_{2}&y_{2}\\
a_{2}&b_{2}
\end{array}
\right]\Big)=
\pi_{\sigma,\rho }\left[\begin{array}{cc}
x_{1}&y_{1}\\
a_{1}&b_{1}
\end{array}
\right]\pi_{\sigma,\rho }\left[\begin{array}{cc}
x_{2}&y_{2}\\
a_{2}&b_{2}
\end{array}
\right]$
can be shown easily using the fact that $\pi_{\sigma,\rho}$ is a ring
homomorphism from $Z_{6}$ to $Z_{3}$.

Note that $\pi_{\sigma,\rho }(\bar{5})=\bar{2}\in U(Z_{3})$, so
$\pi_{\sigma,\rho }$ is an epimorphism from $G_{\sigma}$ to
$G_{\rho}$.

Similarly, the ring homomorphisms $\pi_{\sigma,\beta}:
Z_{6}\longrightarrow Z_{2}$ and $\pi_{\gamma,\beta}:
Z_{4}\longrightarrow Z_{2}$ can be defined respectively satisfying
$\pi_{\sigma,\beta}(\bar{0})=\bar{0}$, $\pi_{\sigma,\beta}(\bar{1})=\bar{1}$,
$\pi_{\sigma,\beta}(\bar{2})=\bar{0}$,
$\pi_{\sigma,\beta}(\bar{3})=\bar{1}$, $\pi_{\sigma,\beta}(\bar{4})=\bar{0}$,
$\pi_{\sigma,\beta}(\bar{5})=\bar{1}$ and
$\pi_{\gamma,\beta}(\bar{0})=\bar{0}$, $\pi_{\gamma,\beta}(\bar{1})=\bar{1}$,
$\pi_{\gamma,\beta}(\bar{2})=\bar{0}$,
$\pi_{\gamma,\beta}(\bar{3})=\bar{1}$. Moreover, the group homomorphisms
$\pi_{\sigma,\beta}: G_{\sigma}\longrightarrow G_{\beta}$ and
$\pi_{\gamma,\beta}:
G_{\gamma}\longrightarrow G_{\beta}$ can be constructed in a similar way.

Finally, we define $\pi_{\beta,\alpha}:  G_{\beta}\longrightarrow
G_{\alpha}$, $\pi_{\rho,\alpha}: G_{\rho}\longrightarrow
G_{\alpha}$, $\pi_{\delta,\sigma}: G_{\delta}\longrightarrow
G_{\sigma}$, $\pi_{\delta,\gamma}: G_{\delta}\longrightarrow
G_{\gamma}$ as the trivial group homomorphisms. Then one has the
following diagram:
\[
\begin{diagram}
\node{G_{\delta}}\arrow{s,r}{\pi_{\delta,\gamma}}
 \arrow{e,t}{\pi_{\delta,\sigma}}\node[1]{G_{\sigma}
}\arrow{s,r}{\pi_{\sigma,\beta}}
\arrow{e,t}{\pi_{\sigma,\rho}}\node[1]{G_{\rho}}
\arrow{s,r}{\pi_{\rho,\alpha}}\\
\node{G_{\gamma}}\arrow{e,t}{\pi_{\gamma,\beta}}
\node[1]{G_{\beta}}\arrow{e,t}{\pi_{\beta,\alpha}}\node[1]{G_{\alpha}}
\end{diagram}
\]

Now, we introduce the multiplication``$\cdot$" in $S$ by
$XW=\pi_{u,uv}(X)\pi_{v,uv}(W)$ if $X\in G_{u}$ and $W\in G_{v}$
for $u,v\in Y$. Then, with this multiplication, $S=\cup_{u\in
Y}G_{u}$ becomes a Clifford monoid related to the semilattice $Y$,
and the only element $e_{\delta}$ of $G_{\delta}$ is the identity
of $S$.

Obviously, $S$ is a finite and noncommutative Clifford monoid.
Then for the cocommutative weak Hopf algebra $kS$ we may obtain
the quantum double $D(S)$ and its quasi-R-matrix $R$ by using the
result in \cite{L1}. We have the decomposition of linear spaces as
follows:
\[
D(S)=(kS)^{op\ast}\infty(kS)=(\oplus_{u\in Y}
kG_{u})^{op\ast}\infty (\oplus_{u\in Y}kG_{u})=\oplus_{u,v\in Y}
((kG_{u})^{op\ast}\infty(kG_{v}))
\]
where $(kG_{u})^{op\ast}\infty(kG_{v})$ means a direct summand of
$D(S)$ and $\infty$ is same in $D(S)$ since
 $(kS)^{op\ast}=\oplus_{u\in Y}(kG_{u})^{op\ast}$ and 
 $kS=\oplus_{u\in Y}kG_{u}$ such that
 for each $u\in Y$, $(kG_{u})^{op\ast}$ is embedded into $(kS)^{op\ast}$ and
 $kG_{u}$ is embedded into $kS$.

 Any $\varphi\in (kG_{u})^{op\ast}$ can be expanded to
 $\overline{\varphi}\in(kS)^{op\ast}$ satisfying
 $\overline{\varphi}(U+V)=\varphi(U)$ for any element $U+V$ of
$kS=\oplus_{v\in Y}kG_{v}$ where $U\in kG_{u}$ and
$V\in\oplus_{v\neq u}kG_{v}$.

Let $u_{1}$, $u_{2}$, $v_{1}$, $v_{2}\in Y$, $X\in G_{v_{1}}$,
$W\in G_{v_{2}}$, $A\in G_{u_{1}}$, $B\in G_{u_{2}}$. Then their
dual elements $\phi_{A}$ and $\phi_{B}$ of $A$ and $B$ are in
$(kG_{u_{1}})^{op\ast}$
 and $(kG_{u_{2}})^{op\ast}$, respectively, where
$\phi_{A}(C)=\left\{\begin{array}{ll}
0  &  \mbox{if $C\in G_{u_{1}}, C\neq A$}\\
1  &  \mbox{if $C=A$}
\end{array}
\right. $
and $\phi_{B}$ is given similarly.

The multiplication of $D(S)$ can be presented  by
\[
(\phi_{A}\infty X)(\phi_{B}\infty W)=(\overline{\phi_{A}}\infty X)
(\overline{\phi_{B}}\infty W)
=\overline{\phi_{A}}\>\overline{\phi_{B}}(X^{-1}?X)\infty XW,
\]
where
\[
\overline{\phi_{A}}\>\overline{\phi_{B}}(X^{-1}?X)=
\left\{\begin{array}{ll}
0  &  \mbox{if $X^{-1}AX\neq B$}\\
\overline{\phi_{A}}=\phi_{A}  &  \mbox{if $X^{-1}AX=B.$}
\end{array}
\right.
\]
By \cite{L1}, the quasi-R-matrix of $D(S)$ is
\[
R=\sum_{s\in S}(1\infty s)\otimes_{k}(\overline{\phi_{s}}\infty 1)
=\sum_{u\in Y}\sum_{g_{u}\in G_{u}}(1\infty g_{u})\otimes_{k}
(\phi_{g_{u}}\infty 1)\in D(S)\otimes D(S).
\]

>From $U(Z_{6})=\{\bar{1}, \bar{5}\}$, $U(Z_{4})=\{\bar{1},
 \bar{3}\}$, $U(Z_{3})=\{\bar{1}, \bar{2}\}$, $U(Z_{2})=\{\bar{1}\}$ and the
fact that $X= \left[\begin{array}{cc}
x&y\\
a&b
\end{array}
\right] \in U(Z_{n}^{2\times2})$ if and only if det$X=xb-ay\in
U(Z_{n})$, it is easy to compute $|G_{u}|$ for each $u\in Y$. We
have $|G_{\delta}|=1, |G_{\alpha}|=1, |G_{\beta}|=6,
|G_{\gamma}|=96, |G_{\rho}|=48, |G_{\sigma}|=288$. It follows that
the number of monomials of $R$ of $D(S)$ is
$|S|=|G_{\delta}|+|G_{\alpha}|+|G_{\beta}|+|G_{\gamma}|+|G_{\rho}|
+|G_{\sigma}|=440$.
Therefore, from the Clifford monoid $S=\cup_{u\in Y}G_{u}$ we have
constructed an example of the quantum doubles of cocommutative
weak Hopf algebras in \cite{L1}.

In the following we  give an example of the quantum doubles of perfect
(noncocommutative) weak Hopf algebras.

For the Clifford monoid $S$ above, $H=kS\otimes(kS)^{\ast}$ is a
finite dimensional non-commutative and non-cocommutative biperfect
weak Hopf algebra with invertible weak antipode $T=T_{S}\otimes
T_{S}^{\ast}$ satisfying $T_{S}(X)=X^{-1}$ and $T^{\ast}
_{S}(f)(X)=f(T_{S}(X))$  for any matrix $X\in S$ and
$f\in(kS)^{\ast}$. The dimension of $H$ is
dim$H=$dim$(kS\otimes(kS)^{\ast})=|S|^{2}=193600$. The quantum
double is given by
\begin{eqnarray*}
D(H)&=&H^{op\ast}\infty H=(kS\otimes(kS)^{\ast})^{op\ast}
\infty(kS\otimes(kS)^{\ast})\\
&=&((kS)^{op\ast}\otimes(kS)^{cop})
\infty(kS\otimes(kS)^{\ast})
=((kS)^{op\ast}\otimes kS)\infty(kS\otimes(kS)^{\ast})\\
&=&\sum_{u,v,p,q\in Y}((kG_{u})^{op\ast}\otimes kG_{v})
\infty(kG_{p}\otimes(kG_{q})^{\ast}),
\end{eqnarray*}
where $kG_{p}\otimes(kG_{q})^{\ast}$ and $(kG_{u})^{op\ast}\otimes
kG_{v}$ are as subspaces of $kS\otimes(kS)^{\ast}$ and
$(kS)^{op\ast}\otimes kS$. The multiplication of $D(S)$ obeys the
formula in Proposition \ref{prop3.1}.

The quasi-R-matrix of $D(H)$ is
$$
R=\sum_{p,\>q\in Y}\sum_{g_{p}\in G_{p},\> g_{q}\in G_{q}}
((1\otimes 1)\infty(g_{p}\otimes\phi_{g_{q}}))\otimes
((\phi_{g_{p}}\otimes g_{q})\infty(1\otimes 1)),
$$
whose number of monomials is $|S|^{2}=193600$.

\vskip.2in
\noindent {\bf Acknowledgement:} We would like to express our sincere
thanks to the anonymous referee for detailed and valuable comments and
suggestions which bring the paper to its present form.

This work is financially supported by
Australian Research Council. Fang Li is also supported by the NSF of
Zhejiang Province of China (Project No.102028) and a
Raybould Fellowship from the University of Queensland.

\end{document}